\documentclass[journal]{IEEEtran}
\usepackage{multirow}
\usepackage{multicol}
\usepackage{lipsum}
\usepackage{graphicx}
\usepackage{graphics}
\usepackage{amsmath}
\usepackage{amsbsy}
\usepackage{blindtext}
\usepackage{tcolorbox}
\usepackage{amssymb}
\usepackage{scalerel}
\usepackage{mathabx}
\usepackage{verbatim}
\usepackage{booktabs}
\usepackage{rotating,tabularx}
\usepackage{mathrsfs} 
\usepackage{hyperref}

\usepackage{algorithm}
\usepackage{algpseudocode}

\usepackage{adjustbox}
\usepackage{soul}
\usepackage{xcolor}
\usepackage{cite}
\usepackage{tikz}
\usepackage{color, colortbl}
\usepackage[first=0,last=9]{lcg}
\definecolor{LightGray}{rgb}{0.7,0.7,0.7}

\usepackage{algorithm}
\usepackage{algpseudocode}
\usepackage{float} 

\usepackage[wby]{callouts}
\usetikzlibrary{shapes.multipart}

\usepackage{array}
\usepackage{makecell}

\allowdisplaybreaks

\usepackage{amsthm}
\theoremstyle{definition}

\theoremstyle{remark}

\usepackage[utf8]{inputenc}
\usepackage[english]{babel}

\usepackage[caption=false,font=footnotesize]{subfig}
\usepackage[T1]{fontenc}
\usepackage{scalerel,stackengine}
\newcommand\reallywidecheck[1]{%
\savestack{\tmpbox}{\stretchto{%
  \scaleto{%
    \scalerel*[\widthof{\ensuremath{#1}}]{\kern-.6pt\bigwedge\kern-.6pt}%
    {\rule[-\textheight/2]{1ex}{\textheight}}
  }{\textheight}%
}{0.5ex}}%
\stackon[1pt]{#1}{\scalebox{-1}{\tmpbox}}%
}
\stackMath
\IEEEoverridecommandlockouts
\IEEEoverridecommandlockouts

\renewcommand{\baselinestretch}{1}

\newcommand*{\rn}{\textcolor{black}}
\newcommand*{\mri}{\textcolor{black}}

\newif\ifarxiv
\arxivtrue

\begin{document}

\title{\LARGE\bf

Assessing EV Charging Impacts on Power Distribution Systems: A Unified Co-Simulation Framework
}

\author{Mohammadreza Iranpour$^{\ast}$, Mohammad Rasoul Narimani$^{\ast}$, Xudong Jia$^{\dagger}$
\thanks{${\ast}$: Department of Electrical and Computer Engineering, California State University Northridge (CSUN).
mohammadreza.iranpour.443@my.csun.edu,
Rasoul.narimani@csun.edu, 
} 
\thanks{${\dagger}$: Department of Civil Engineering and Construction Management, California State University Northridge (CSUN).
xudong.jia@csun.edu.
This work was supported by University of California Office ofthe President Climate Action under grant number \#R02CP6948.}
}

\maketitle

\begin{abstract}

\rn{The growing adoption of electric vehicles (EVs) is expected to significantly increase the load on electric power distribution systems, many of which are already operating near their capacity limits. To effectively address this challenge, this paper presents a comprehensive framework for analyzing the impact of large-scale EV integration on power distribution networks. The proposed framework utilizes the open-source simulator OpenDSS to build a detailed, scalable model of electric distribution systems. It leverages high-fidelity synthetic data from the SMART-DS project, modeling three feeders from an urban substation in San Francisco, California, down to the household level. A key contribution of this work is its ability to identify critical components within the distribution system that may require upgrades due to increased EV-related loads. The framework incorporates advanced geospatial visualization using QGIS, which enhances understanding of how additional EV charging demands affect different parts of the grid, enabling stakeholders to pinpoint where infrastructure reinforcements are most needed. To ensure realism in modeling EV loads, the framework integrates projected load profiles derived from U.S. Department of Energy data, which accounts for multiple influencing factors such as vehicle types, usage patterns, charging behaviors, and adoption rates. Furthermore, the use of realistic, large-scale synthetic data ensures the model’s applicability to real-world planning and decision-making. The framework supports various simulation scenarios, including light to heavy charging loads and comparisons between distributed and centralized charging patterns, offering a practical tool for utilities and planners to prepare for widespread EV adoption. In addition, the modular structure of the framework allows for easy adaptation to different geographic locations, feeder configurations, and EV adoption rates. This flexibility makes it well-suited for future studies on evolving grid conditions.}

\end{abstract}


\section{Introduction}
\label{Introduction}

\rn{Declining battery costs and federal incentives have significantly accelerated the adoption of electric vehicles (EVs) in recent years, as reported by the U.S. Department of Energy. Consequently, a growing number of EV owners are charging their vehicles at home, leading to a noticeable rise in residential charging station installations, according to data from the U.S. Energy Information Administration. This shift marks a fundamental change in electricity consumption patterns, especially at the distribution level, where power is delivered directly to consumers. As more EVs are plugged in during evening hours, already a time of peak residential electricity demand, the risk of overloading local distribution networks increases. This can result in operational issues such as voltage deviations, transformer overloading, and increased wear on system components. If not properly managed, these challenges could compromise the reliability and efficiency of the grid. As such, understanding and mitigating the impacts of widespread home EV charging is becoming increasingly important for utilities, policymakers, and grid operators~\cite{ahmad2022optimal}.}

\rn{Power distribution systems have limited capacity to handle overloads, so the growing demand from EV charging may require upgrading existing infrastructure or adding new capacity~\cite{Mogos2021}. Charging stations vary in charging speed, such as slow and fast chargers, which place different levels of demand on the grid and thus have varying impacts on power system~\cite{Hu2012}. A detailed overview of how EV charging stations affect the power grid is provided in~\cite{rahman2022comprehensive}. These impacts can influence grid stability, voltage profiles, and peak load management, making strategic planning for EV integration essential.}

\rn{Extensive research has focused on mitigating the impact of EV charging on the power grid. One proposed solution, particularly for fast charging stations, is to connect them to medium-voltage transmission lines. This approach takes advantage of higher voltage levels, which can accommodate additional loads with minimal effect on line characteristics, since the power demand from EV chargers is relatively small compared to the total transmission system load~\cite{Li2019}. Another widely studied strategy is smart charging, which dynamically schedules EV charging based on grid conditions, electricity prices, and the availability of renewable energy. The work in~\cite{crozier2020opportunity} presents a framework that considers vehicle usage patterns, energy consumption, and grid topology to simulate various charging scenarios using conditional probability models. The goal is to facilitate a smooth transition to large-scale EV adoption without overloading the grid. EV charging management strategies can generally be categorized into centralized and decentralized approaches. In centralized control, a single entity such as a utility or aggregator coordinates the charging schedule for multiple EVs. In contrast, decentralized control allows individual EVs to independently schedule charging based on local information, such as real-time electricity prices~\cite{CHENG2018175}. Both approaches aim to reduce peak demand~\cite{ramadan2018smart, narimani2017energy, narimani2017}, minimize energy losses~\cite{nafisi2015two, narimani2015effect}, maintain phase balance in the distribution network~\cite{weckx2015load}, and increase the integration of renewable energy sources~\cite{clairand2017tariff}. Most charging algorithms are designed to meet EV charging requirements without exceeding network constraints, while some incorporate additional factors such as battery temperature or grid performance metrics to further optimize charging efficiency and system reliability~\cite{nour2019smart}. While prior studies have explored strategies for managing EV charging impacts and proposed various control methods and modeling approaches, most focus on either limited-scale networks or simplified load assumptions.}

\rn{Several studies employ high-level aggregated EV load models to rapidly screen distribution-system impacts without performing full time-series network simulations. One approach formulates EV charging demand as nodal injections in a quasi-static power-flow framework and optimizes charging schedules for low-voltage feeders \cite{richardson2011optimal}. Another method represents an entire EV fleet as a stochastic demand block based on arrival and departure statistics, enabling fast Monte Carlo assessment of thermal and voltage violations \cite{paris2016aggregate}. A further study develops a probabilistic aggregated model that captures the variability of charging events for scenario-based planning applications \cite{yu2016stochastic}.} However, The rapid growth in \mri{EVs} necessitates analyses that utilize high-fidelity models of electric distribution systems. High-fidelity simulations can accurately assess the impacts of \mri{EVs} on power grids\mri{~\cite{Wang2020}}. 
\rn{Multiple open-source and commercial platforms have been leveraged to evaluate the impacts of EV charging on distribution networks. An open-source GridLAB-D simulator is used to model smartphone-controlled EV charging and assess its effects on feeder voltages and losses~\cite{zhan2017modeling}.  The Python library pandapower, built atop PYPOWER, has also been applied to low-voltage network studies including aggregated EV charging scenarios \cite{thurner2018pandapower}.  Likewise, the PyPSA framework enables high-resolution time-series power-flow simulations with stochastic EV loads \cite{brown2018pypsa}.  On the commercial side, DIgSILENT PowerFactory to quantify thermal and voltage violations in fast-charging clusters has been employed in \cite{khorasani2020powerfactory}}.

\rn{Building on prior efforts, this paper presents a fully integrated framework that overcomes key limitations in existing approaches to EV impact analysis. While previous studies have provided useful approximations, they often rely on simplified EV load models, lack spatial resolution, and do not support integration with standard load profiles (SLPs) or high-resolution smart meter data. Furthermore, most open-source and commercial tools do not offer native support for combining EV-specific load patterns with existing residential and commercial load data, limiting their applicability for detailed, realistic grid impact assessments.
To address these gaps, we develop a high-fidelity simulation framework that combines OpenDSS for detailed distribution system modeling, Julia for flexible and high-performance scripting, and GIS-based tools for spatially-aware EV station placement and load characterization. This integrated approach enables seamless, reproducible, and geographically-resolved analysis of EV charging impacts. Our framework captures the combined effects of EV and non-EV loads, evaluates line overloads and voltage deviations across varying EV penetration scenarios, and provides actionable insights for utility operators and grid planners preparing for large-scale EV adoption.}



\rn{This paper presents a high-fidelity model for simulating EV charging across multiple feeders in power distribution systems, using the OpenDSS simulation platform~\cite{dugan2011open}. By modeling distribution networks down to the household level, the framework offers a level of detail that exceeds what is typically achievable through analytical methods. This granularity enables the model to serve as a validation tool for evaluating modeling assumptions and testing control strategies developed using simplified grid representations. Recent advances in synthetic power system datasets have also facilitated detailed studies of transmission and distribution networks~\cite{li2020building}. Our work builds on one such dataset from the SMART-DS project~\cite{krishnan2017smart}, modeling numerous feeders from an urban substation in San Francisco.}

\rn{A key contribution of this work is its ability to identify critical components within the distribution network that may require reinforcement under increased EV charging demand. The framework incorporates visualization tools to support intuitive analysis of overload and voltage issues, enabling planners and utilities to pinpoint infrastructure vulnerabilities. Additionally, the model supports various simulation scenarios, including different charging intensities. Its modular design allows easy adaptation to different geographic regions, feeder configurations, and EV adoption levels, making it a practical and flexible tool for both current analysis and future research on evolving grid conditions.}

\rn{The remainder of this paper is organized as follows. Section~\ref{sec:Methodological Framework} describes the structure of the modeling framework and the sources of data used for its components. Section~\ref{sec:model} details the proposed high-fidelity model for analyzing the impact of EV integration on low-voltage distribution systems. Simulation results are presented and discussed in Section~\ref{sec:results}. Finally, Section~\ref{sec:conclusion} summarizes the main findings and outlines directions for future work.}

\section{Methodological Framework}
\label{sec:Methodological Framework}

\rn{In this section, we describe the key components used in the proposed methodology for assessing the impact of EV loads on power distribution systems.}

\subsection{Incorporating EV loads into Load Profile}
\label{sec:load_profile}

\rn{Assessing current electricity usage in the power distribution system is essential for understanding the potential impact of EV charging. To support this analysis, we use data from the SMART-DS dataset~\cite{krishnan2017smart}, focusing on a sub-region within the extended San Francisco Bay Area model. The SMART-DS dataset includes both peak planning loads and time-series load data. Peak loads, provided in real (kW) and reactive (kVar) power, represent the maximum expected demand on the network. The time-series data consist of 15-minute resolution load profiles over an entire year, capturing fine-grained variations in demand throughout different times and seasons. SMART-DS includes load profiles for residential, commercial, and industrial consumers~\cite{krishnan2017smart}.
In this study, we specifically use the residential load profiles to evaluate how EV charging may affect low-voltage (LV) distribution systems. Each time-series profile in SMART-DS also includes a breakdown of end-use components~\cite{krishnan2017smart}, which helps in analyzing the contribution of specific load types, such as HVAC, lighting, and appliances. This high-fidelity data enables a more accurate and realistic assessment of baseline load behavior prior to the introduction of EV charging.}
\rn{A key advantage of using the SMART-DS dataset is its grounding in realistic, utility-grade data that reflect actual customer usage patterns and grid configurations. The dataset is designed to support distribution system analysis by providing high-resolution, geographically detailed, and temporally accurate information. This makes it especially useful for simulating the dynamic interactions between traditional loads and emerging technologies like EVs. For our analysis, we extract data for several feeders within the "PU15" sub-region, a small urban area located in San Francisco, from the extended Bay Area model. This area provides a representative environment for studying the potential stress that widespread EV adoption could place on urban distribution infrastructure.}

\subsection{
\mri{Modeling and Characterization of EV Charging Loads}}
\label{sec:Charging Characteristic}


\rn{For assessing the impacts of EV charging on distribution systems it is essential to consider both the performance characteristics and the geospatial distribution of charging stations. We draw on the U.S. Department of Energy’s Alternative Fuels Data Center (AFDC) \cite{afdc2025}, launched in 1991 in response to the Alternative Motor Fuels Act and Clean Air Act Amendments and continuously reviewed by subject-matter experts, as our primary source of data related to the EV station behavior and load profile. To evaluate how charging stations perform, we use the AFDC’s EVI-X Toolbox. This tool helps us account for key factors such as the number of electric vehicles, how far they travel each day, whether drivers can charge at home or at work, their charging habits, and the types of charging stations available. For spatial allocation we employ the AFDC’s Alternative Fueling Station Locator, which displays publicly accessible stations across the United States and Canada, offers advanced filters by fuel type, charger power output (50-350 kW), charger level (Level 1, Level 2, DC FAST Chargers), station status (including temporarily unavailable), and relies on station data gathered and verified by the National Renewable Energy Laboratory (NREL) through trade media, Clean Cities coalitions, provider submissions, original equipment manufacturers, and industry collaboration.
In the following sections, we present a detailed description of the approaches and assumptions underlying these factors.}

\subsection{\rn{Power Distribution Network and Power Flow Analysis}}
\label{sec:Power Distribution Network}

\rn{To assess the impact of EV charging on power distribution systems, we use detailed data from the SMART-DS dataset. This dataset allows us to perform power flow analyses across various levels of the grid, including entire regions, individual substations, and specific feeders. For feeder-level studies, we simulate a single feeder with a fixed voltage at the feeder head. To represent EV charging behavior, we adjust household loads based on typical charging patterns associated with Level 2 chargers. By scaling residential loads to reflect these patterns, we can simulate how EV charging affects the performance of the distribution network.}

\rn{We perform power flow simulations using OpenDSS~\cite{OpenDSS}, a widely used simulation tool for electric power distribution systems. These simulations help calculate bus voltages and network losses, allowing us to evaluate how EV charging influences grid stability and efficiency. Our simulation model is also flexible and can be easily customized to test different loading conditions across various time windows. To enhance computational efficiency and integration with custom analysis tools, we use OpenDSSDirect, a Python-based interface for OpenDSS. OpenDSSDirect enables direct control of active and reactive loads, offering greater flexibility in modeling a wide range of operating scenarios. One key advantage is that it allows us to read time-series data, such as load shapes, directly from parquet files, eliminating the need to load all data into memory. This significantly speeds up power flow computations, especially for large-scale simulations. By combining realistic load data with efficient simulation tools, our methodology provides a robust framework for analyzing the impacts of EV charging on distribution networks, supporting better planning and management of grid infrastructure.}

\section{High-level Fidelity Model}
\label{sec:model}

\begin{figure*}
    \centering
\includegraphics[scale=0.33,trim=4cm 0.3cm 4cm 0.5cm,clip]{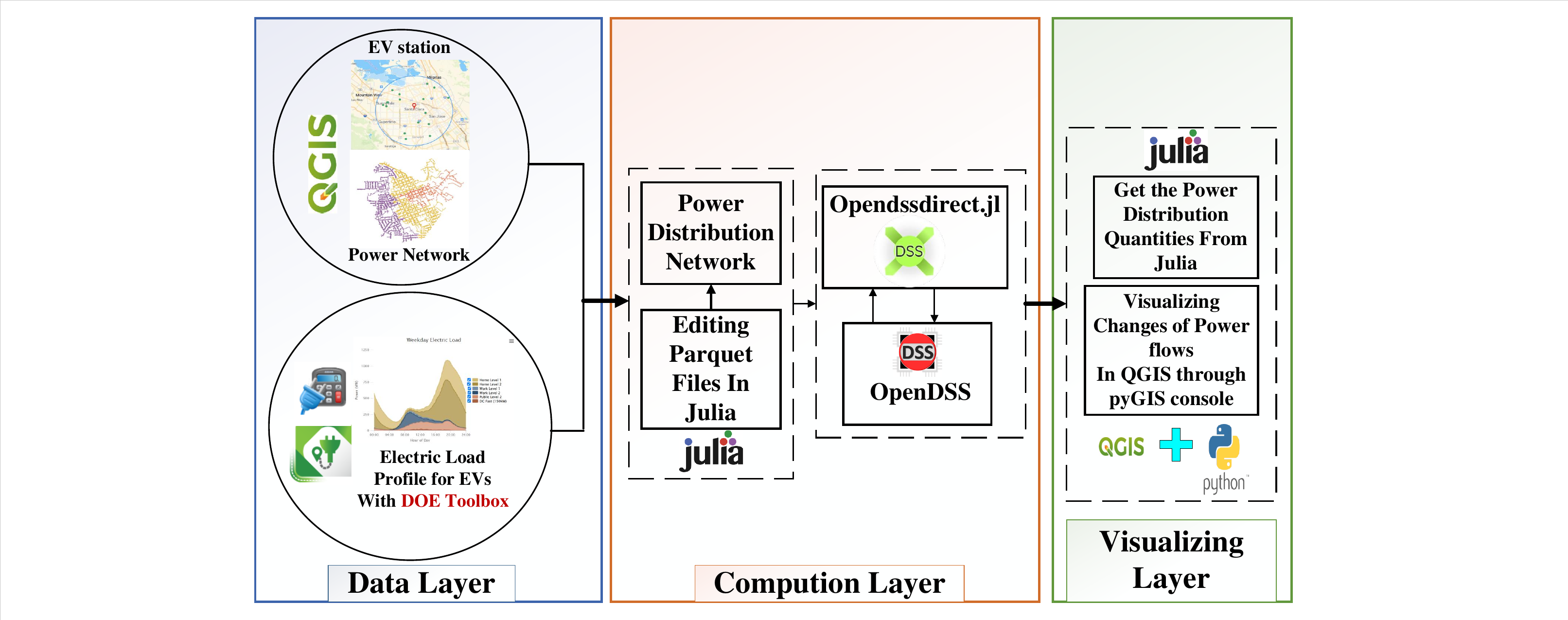}
    \caption{Schematic of the proposed framework for analyzing the impact of EV chargers on power distribution systems. The framework is organized into three integrated layers, i.e., data, computation, and Visualization. Geographic and electrical data are first collected and processed using GIS and DOE datasets. EV load profiles are then estimated using the EVI-X modeling suite and assigned to load buses. Finally, time-series power flow simulations are performed using OpenDSSDirect to evaluate system-level impacts, including voltage profiles, line loading, and energy losses.}
    \label{fig:procedure}
\end{figure*}

\rn{In this study, we propose a multi-layered methodology to assess the impact of EV chargers on power distribution systems. The approach is structured into three main layers: the data layer, the computational layer, and the visualization layer. Each layer plays a distinct role and utilizes specific tools to support the overall analysis. A detailed description of these layers is provided in the following sections. The general schematic of the proposed methodology is shown in Fig.~\ref{fig:procedure}.}

\subsection{Data Layer}
\label{sec:data layer}

\rn{In our framework for analyzing the impact of EV chargers on power distribution systems, multiple data layers are integrated to accurately evaluate EV load effects on distribution system load buses. Geographic and electrical characteristics of the distribution network are obtained using GIS-based data, specifically from JSON files provided by the SMART-DS dataset (see Section~\ref{sec:load_profile}). Additionally, we incorporate EV charging station locations from the U.S. Department of Energy (DOE) database \cite{AFDC_EV_Stations}, which includes stations with power ratings ranging from 50 kW to over 350 kW. These datasets are imported into QGIS using CSV files for spatial analysis and visualization.}

\rn{To estimate the electrical load demand of these stations, we use the Electric Vehicle Infrastructure (EVI-X) Toolbox developed by the NREL \cite{DOE_EV_Toolbox}, as described in Subsection~\ref{sec:Charging Characteristic}. EVI-X supports planning and modeling of EV infrastructure at various scales, from national to site-specific levels. The EVI-Pro Lite tool within this suite estimates the charging infrastructure requirements for a given region based on factors such as average daily vehicle mileage, EV adoption rates, home charging availability, and user charging behavior.}
\rn{Using these inputs, we apply a uniform distribution method to allocate the estimated EV load across the charging stations, adjusting for their respective power capacities. These loads are then mapped to the nearest load buses defined in the SMART-DS JSON files. This enables a detailed assessment of the impact of EV charging on the power distribution system.}

\subsection{Computational Layer}
\label{sec:computational layer}

\rn{As outlined in the previous subsection, the first step in performing power flow analysis on the distribution system is assigning the EV loads to the appropriate load buses. To achieve this, we utilize the OpenDSSDirect package, a flexible and open-source tool specifically developed for the analysis of power distribution systems. OpenDSSDirect~\cite{OpenDSSDirect} serves as a high-level interface to the OpenDSS simulation engine via the DSS C-API and supports both Python and Julia environments. In this study, we adopt the Julia implementation (OpenDSSDirect.jl), which allows for efficient integration into custom analysis workflows and offers robust cross-platform compatibility.}

\rn{To simulate the temporal evolution of EV charging demand, we developed a quasi-static time-series simulation module that performs power flow calculations over a full annual cycle, consisting of 8,760 hourly time steps. The simulations are executed on a distribution feeder model represented in OpenDSS format. This requires two key components: a detailed feeder model describing the network topology and electrical parameters, and a corresponding time-series load profile capturing the variation in demand throughout the year. These elements are managed through a set of DSS files, with the simulation coordinated by a master DSS script executed through the OpenDSSDirect engine.}

\rn{OpenDSSDirect supports direct assignment of real and reactive power values (kW and kVAR) for each load, which eliminates the need to pre-load external load shape files into memory. This feature enables faster execution of time-series simulations. For this study, the EV load profiles are obtained from Parquet files located in the load-data directory. These profiles represent the hourly charging demand at each station and are manually annotated in the corresponding DSS load definitions within the opendss-no-loadshapes directory.}

\rn{This modeling approach allows us to capture the dynamic behavior of EV loads with high temporal resolution and directly observe their impacts on the power distribution system. Specifically, we can assess how different patterns of EV charging affect key performance indicators such as voltage magnitude deviations, line loading, transformer utilization, and overall system losses. The ability to modify load profiles and analyze system response in an integrated environment also enables scenario-based studies, such as evaluating the effects of increased EV adoption or different charging strategies (e.g., time-of-use rates or managed charging).}
\rn{The full procedure for assigning EV loads to the nearest load buses and executing power flow simulations using OpenDSSDirect.jl is detailed in Algorithm~\ref{alg:load_assignment}. This capability forms a foundational layer of our framework for evaluating the operational impacts of EV integration into existing distribution networks and supports further studies on system planning, resilience, and optimization under high EV penetration scenarios.}

\begin{algorithm}
\caption{Process of Assigning Loads to Nearest Buses and Running Power Flow}\label{alg:load_assignment}
\begin{algorithmic}[1]
\State \textbf{Input:} 
\State EV station locations (CSV files), Load bus locations (JSON files), Load profile for each EV station, Master DSS files.
\State \textbf{Output:} Updated load assignments for buses, Power flow results.

\State \textbf{Step 1:} \textit{Load EV station data from CSV files.}
\State \quad \textbf{For each CSV file:}
\State \quad \quad Read station name, latitude, longitude.
\State \quad \quad Store EV stations in a dictionary with the station name as key, and coordinates as value.

\State \textbf{Step 2:} \textit{Load bus coordinates from JSON files.}
\State \quad \textbf{For each JSON file:}
\State \quad \quad Read the bus load names, latitude, and longitude.
\State \quad \quad Store buses in a dictionary with the bus name as key, and coordinates as value.

\State \textbf{Step 3:} \textit{Assign the EV load profile to the nearest bus.}
\State \quad \textbf{For each EV station:}
\State \quad \quad Find the nearest bus based on the minimum distance. 
\State \quad \quad Assign the estimated load profile for each EV station to its corresponding nearest bus.

\State \textbf{Step 4:} \textit{Update Load DSS file.}
\State \quad For each load bus:
\State \quad \quad Update the load value in the "Loads.dss" file based on the assigned EV load.

\State \textbf{Step 5:} \textit{Run Power Flow based on Master DSS file.}
\State \quad \textbf{For each Master DSS file:}
\State \quad \quad Call the Master DSS file using OpenDSSDirect.jl.
\State \quad \quad Solve the power flow using OpenDSSDirect.
\State \quad \quad Collect and store power flow results, including losses and flows on lines.

\State \textbf{Step 6:} \textit{Return Results.}
\State \quad Output the updated load buses with assigned EV loads and power flow results (losses, line flows).
\end{algorithmic}
\end{algorithm}

\subsection{Visualization Layer}
\label{sec:visualizing layer}


\rn{In this section, we describe the visualization layer of our proposed framework, which uses Julia, Python, and QGIS to assess the impact of EV charging loads on GIS-based electrical feeders and lines. Following the power flow analysis performed with the OpenDSSDirect package in Julia, we update and regenerate GeoJSON files to reflect changes in the distribution system resulting from EV integration.}

\rn{Our visualization process consists of several steps. First, we identify high-ampacity lines by filtering those in the GeoJSON dataset that exceed a predefined ampacity threshold. We then compare power flow values before and after incorporating EV station loads by analyzing both the original and modified power flow JSON files. The percentage change in power flow on each high-ampacity line is calculated and used to determine the visual representation of that line. Specifically, lines are categorized by the magnitude of change: gray for changes less than 0.05\%, green for 0.05–50\%, magenta for 50–80\%, and red for changes of 80\% or more. This classification enables a clear and intuitive visual interpretation of where the network experiences the most significant stress due to EV demand.}

\begin{figure}
    \centering 
\includegraphics[scale=.45,trim=0.4cm 5cm 1cm 6cm,clip]{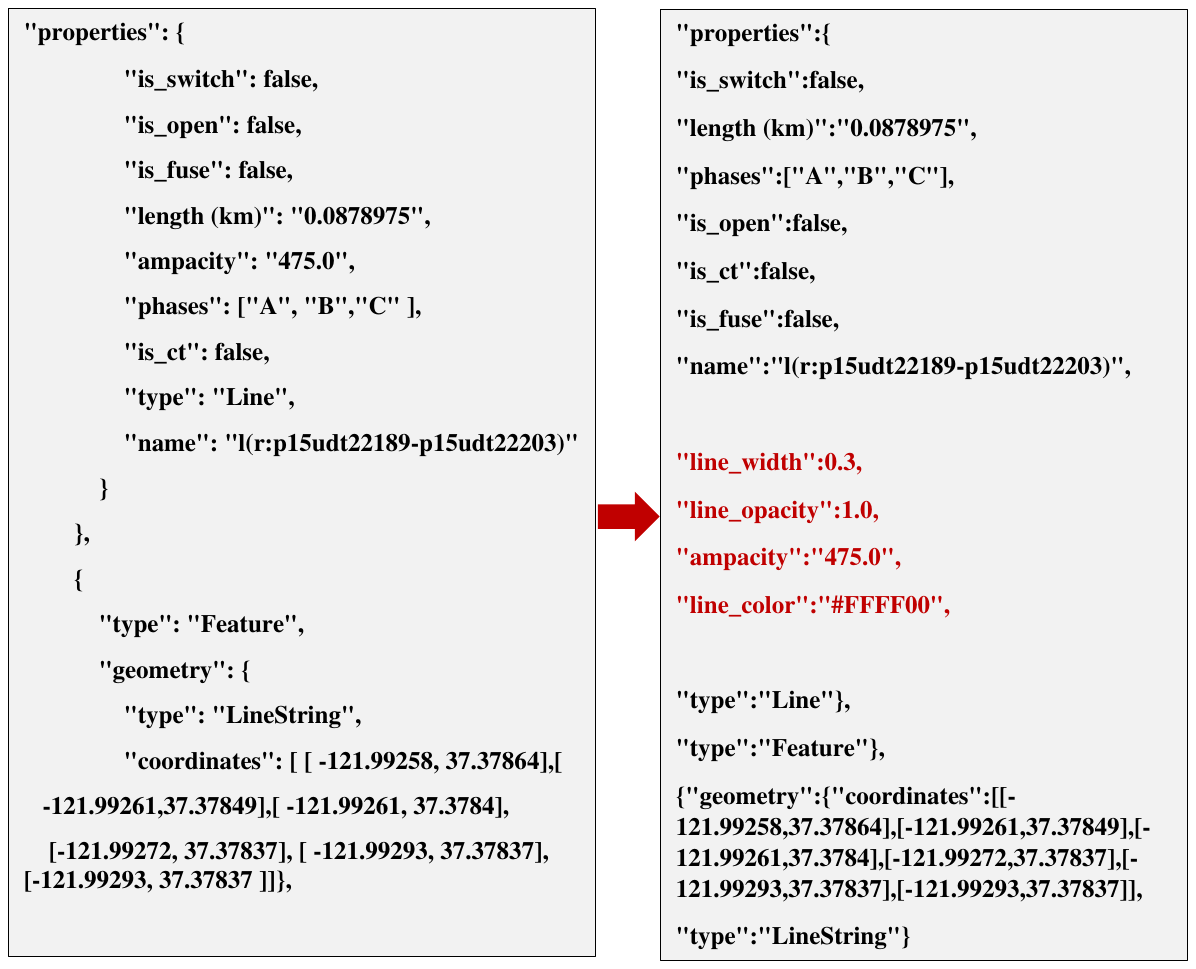}
\centering 
\caption{\rn{Modification of a JSON dictionary entry representing a power distribution line. The original entry (left) includes basic electrical properties, while the updated entry (right) incorporates additional visualization attributes such as ``line\_width'', ``line\_opacity'', and ``line\_color'' to reflect changes in power flow for geospatial rendering.}}
    \label{fig:code}
\end{figure}

\rn{To enhance and customize the visualization, we also add new properties to the GeoJSON files, including line thickness, opacity, and color, which can be directly interpreted by QGIS. These attributes are updated using Python scripts executed in the PyQGIS console. The resulting GIS map allows for interactive exploration of how EV loads affect the distribution system, providing insights into localized stress and potential areas of concern under different charging scenarios.}
\rn{Through this methodology, our framework provides a clear, data-driven approach to assessing the impact of EV charging infrastructure on the distribution network, highlighting potential congestion points and feeder overload risks. A simple example of the JSONs files for before and after changing the loads in the system has been shown in Fig.~\ref{fig:code}. The overall visualization process is outlined in Algorithm\ref{alg:QGIS}.}

\rn{The proposed framework serves as a practical and modular tool for evaluating the impact of large, distributed EV charging loads on power distribution systems. Its layered structure, comprising data acquisition, load modeling, and power flow analysis, enables system operators and planners to systematically assess how various EV-related factors, such as station placement, charging behavior, and load intensity, influence grid performance. This modularity also allows power system researchers to isolate and study the effects of specific variables, supporting independent and targeted analysis. In the following section, the proposed algorithm is applied to two distinct case studies that utilize large-scale, synthetic data to demonstrate the framework's capabilities and effectiveness in capturing the operational impacts of widespread EV integration.}

\section{Results and Analysis}
\label{sec:results}

\rn{In this section, we evaluate the capability of the proposed framework to assess the impact of EV loads on the power grid. To this end, we apply the framework to a range of deterministic scenarios that reflect different EV load patterns. These scenarios are generated using the EVI-X toolbox developed by the U.S. Department of Energy \cite{DOE_EV_Toolbox}. As described in the previous section, the toolbox offers diverse EV charging demand profiles based on urban area characteristics, EV penetration levels, and user charging behaviors. For this study, we use it to derive realistic EV load profiles representative of the San Jose region.}

\begin{figure}
    \centering 
\includegraphics[scale=.32,trim=0.4cm 2.5cm 3.5cm 2cm,clip]{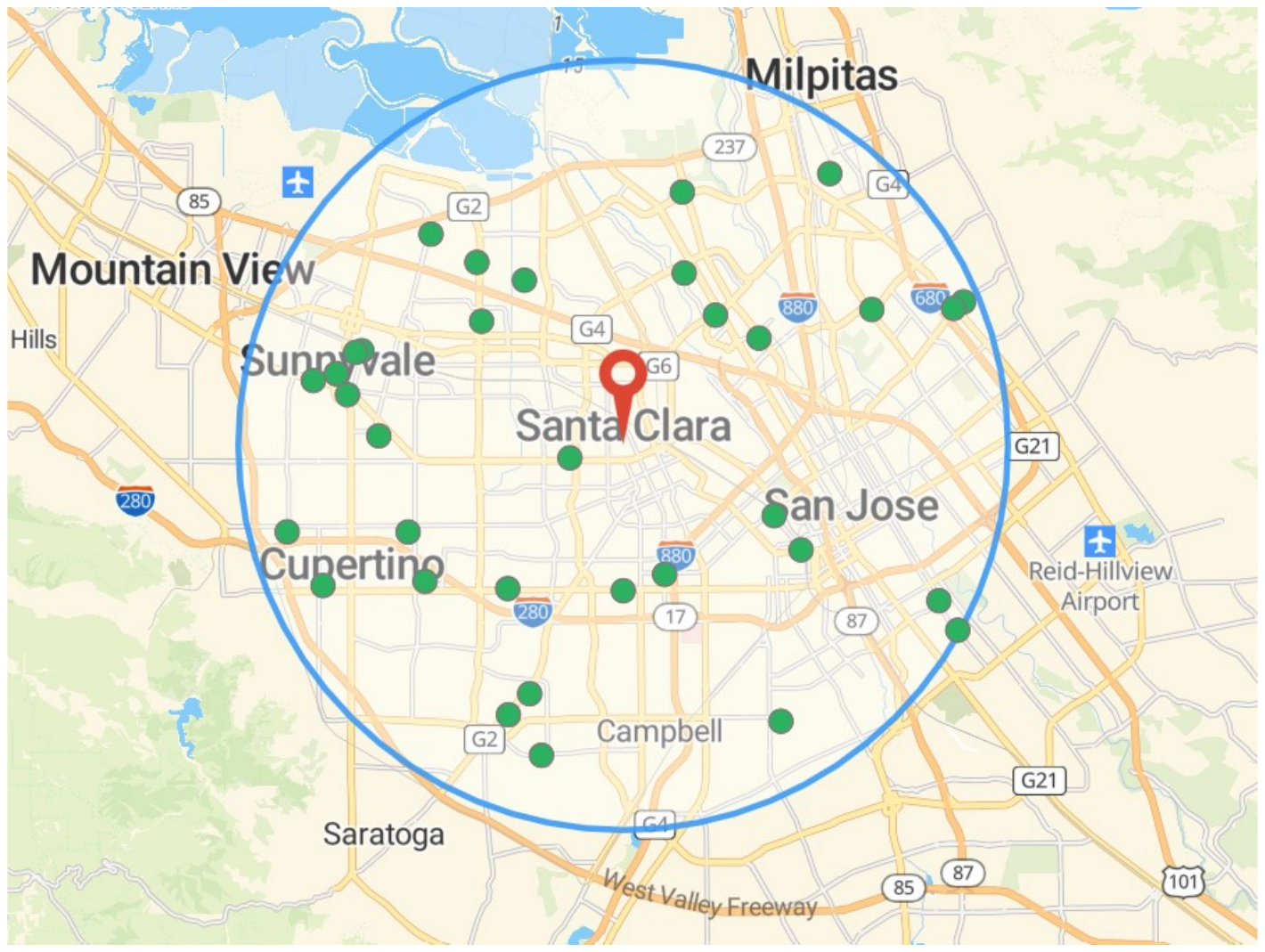}
\centering 
\caption{\rn{Geographic distribution of EV charging stations with power outputs ranging from 50 kW to 150 kW, located within and around San Jose, Sunnyvale, Cupertino, and Santa Clara. The blue circle indicates the area of interest used for spatial analysis.}}
    \label{fig:ev location}
\end{figure}

\rn{To enhance the realism and accuracy of the simulation, we incorporate large-scale, GIS-based data on public EV charging stations from the Alternative Fuels Data Center (AFDC) \cite{AFDC_EV_Stations}. This dataset includes the geographic location, power capacity, and infrastructure specifications of existing EV charging stations. For effective modeling, stations are categorized into four capacity levels: Level 1 (<50 kW), Level 2 (50–150 kW), Level 3 (150–350 kW), and Level 4 (>350 kW). These classifications provide a structured basis for assessing grid impacts according to charging infrastructure.
The total EV charging demand during the peak hour, identified from the 24-hour load profiles generated by the DOE toolbox, is proportionally allocated across 951 public charging stations located in San Jose, Sunnyvale, and Cupertino, based on each station’s capacity. For example, stations with power capacities between 50 kW and 150 kW are illustrated in Fig. \ref{fig:ev location}; note that each point on the map may represent multiple co-located chargers.
By integrating detailed, real-world infrastructure data at scale, the simulation achieves greater fidelity in capturing spatial and operational characteristics of EV charging demand. Moreover, the implementation is designed to be customizable, enabling researchers to investigate how different factors, such as station density, or user behavior, affect the aggregate load profile and the resulting impacts on the power grid.}

\begin{algorithm}
\caption{Visualizing Layer for EV Load Impact on Power Distribution Network}\label{alg:QGIS}
\begin{algorithmic}[1]
\State \textbf{Input:} GeoJSON files of the distribution network
\State \textbf{Input:} Power flow data from OpenDSSDirect (after assigning EV loads)

\State \textbf{Step 1: Load Distribution Network Data}
\State Load the distribution network data from GeoJSON files containing network topology and feeder information in Julia.
\State Parse the network features to extract the lines' coordinates, ampacity, and other necessary attributes.

\State \textbf{Step 2: Compute Power Flow Change and Assign Colors}
\For{Each line in the network}
    \State Retrieve the power flow value for the current line from the OpenDSSDirect results.
    \State Retrieve the ampacity of the current line from the network data.
    \State Calculate the percentage change in power flow compared to its original state.
    \If{Percentage change $<$ 0.05\%}
    \State Assign \textbf{Gray}  color (\#808080)
\ElsIf{0.05\% $\leq$ Percentage change $<$ 10\%}
    \State Assign \textbf{Green} color (\#00FF00)
\ElsIf{10\%  $\leq$ Percentage change $<$ 50\%}
    \State Assign \textbf{Blue}  color (\#0000FF)
\ElsIf{50\%  $\leq$ Percentage change $<$ 80\%}
    \State Assign \textbf{Pink} color (\#FF00FF)
\Else
    \State Assign \textbf{Red}   color (\#e31a1c)
\EndIf
    \State Update the line color property in the network data.
\EndFor

\State \textbf{Step 3: Apply Modified Line Colors to Visualization}
\For{Each line in the GeoJSON layer}
    \State Check if the line has been assigned a new color based on its power flow change.
    \If{Line has a modified color}
        \State Update the line visualization properties (e.g., line color) in QGIS.
    \Else
        \State Apply default line style.
    \EndIf
\EndFor

\State \textbf{Step 4: Visualize the Distribution Network in QGIS}
\State Add the modified network data with updated line colors to the QGIS project.
\State Apply a style rule to color the lines based on their assigned categories.
\State Render the updated map layer in QGIS, showing the impact of EV load on the distribution network.

\end{algorithmic}
\end{algorithm}

\rn{To incorporate EV charging demand in a realistic power distribution system, we use the SMART-DS synthetic electrical network dataset \cite{palmintier_smart-ds_2020}. Developed by the NREL, SMART-DS provides detailed, large-scale synthetic distribution and transmission networks for various U.S. cities, including San Francisco. In this study, we focus on the PU15 region from the San Francisco dataset, which includes the distribution and transmission network topology, time-series load data for the San Jose area, and geospatial JSON files containing GIS coordinates of network elements. These data enable us to model the integration of EV loads into the power system and evaluate their impacts. A one-line diagram of the PU15 zone, showing the proposed EV charging station locations across four deployment levels, is presented in Fig.~\ref{fig:addevtoPU15}. Then, using the GIS-based locations of EV charging stations and power infrastructure elements from the SMART-DS dataset, each EV station is assigned to its nearest transformer in the network based on the algorithm described in Algorithm~\ref{alg:load_assignment}. After establishing the loads, the OpenDSS network files are updated by incorporating the EV load demand at the corresponding transformers using OpenDSSDirect.jl in the Julia environment. }

\rn{Next, to evaluate the impact of EV integration on the power system, we perform a before-and-after power flow analysis using OpenDSSDirect.jl, following the assignments defined by Algorithm~\ref{alg:load_assignment}. OpenDSSDirect.jl offers a high-performance interface for power system simulations, making it well-suited for analyzing large-scale datasets such as SMART-DS. Finally, to visually assess the impact of EV loads on the power system, we modify the color of transmission lines based on their utilization levels, following the procedure outlined in Algorithm~\ref{alg:QGIS}. This visualization aids in understanding how EV integration affects the grid and helps identify areas that may require upgrades. In addition to the geospatial analysis, we perform statistical evaluations to quantify system impacts. Specifically, we generate histograms of power flow distributions before and after EV load integration to illustrate how line flows change due to the added demand. We also compare total system power losses across different scenarios to evaluate the effects of EV integration on overall system efficiency.}

\rn{Before implementing our approach on different scenarios, we first provide insight into key parameters that influence the generation of EV load profiles. Table~\ref{Table:FACTORS} lists several factors that affect EV charging demand, serving as the basis for constructing diverse scenarios. These parameters include the number of plug-in EVs in the fleet, average daily mileage, ambient temperature, vehicle type, charging behavior, and access to charging infrastructure. The assumptions are derived from data available through the DOE EV Toolbox \cite{DOE_EV_Toolbox}. Each factor, along with its possible states and descriptions, supports the development of realistic EV demand curves. This enables the proposed approach to assess the impact of EV loads on the grid under a variety of conditions. In this study, we focus on two representative scenarios derived from these parameters.}

\subsection{First Scenario}
\label{sec:first scenario}

\begin{figure}
    \centering 
\includegraphics[scale=0.45,trim=0.1cm 0.1cm 0.1cm 0.1cm,clip]{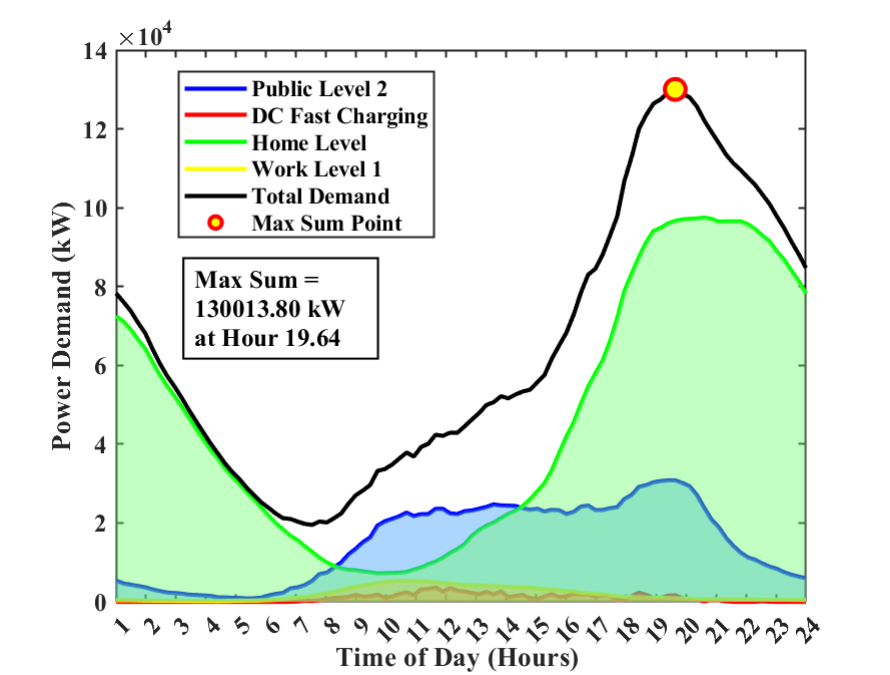}
\centering 
\caption{\rn{24-hour EV charging demand profile for Scenario 1, based on a fleet of 350,000 vehicles in the San Jose–Sunnyvale–Cupertino area. The fleet includes 50\% BEVs and 50\% PHEVs, with most vehicles (80\%) charging at home using a mix of Level 1 and Level 2 chargers. Workplace charging also uses a 50/50 mix of Levels 1 and 2. Home charging starts immediately at low speed, while workplace charging starts immediately at full speed. The peak total demand reaches about 130 MW around hour 19.38.}}
    \label{fig:scenario1}
\end{figure}

\rn{In the first scenario, we apply the following assumptions based on Table~\ref{Table:FACTORS}: The EV fleet consists of 350,000 vehicles. Each vehicle travels an average of 25 miles per day, representing lower daily usage. The average ambient temperature is set to 80°F, reflecting warmer climate conditions. The fleet is composed of 50\% battery electric vehicles (BEVs) and 50\% plug-in hybrid electric vehicles (PHEVs), with 50\% of all plug-in vehicles classified as sedans. For workplace charging, we assume an even mix of 50\% Level 1 and 50\% Level 2 charging infrastructure, balancing slow and faster charging options. All vehicles (100\%) have access to home charging, also with a 50\% Level 1 and 50\% Level 2 distribution. Additionally, 80\% of the fleet primarily charges at home, with workplace charging as a secondary option. The home charging strategy is set to ``Immediate, as slow as possible,'' distributing the load evenly over the parking period. In contrast, the workplace charging strategy is ``Immediate, as fast as possible,'' ensuring vehicles begin charging at maximum speed upon arrival. Based on these assumptions, we generated load demand curves for Level 1, Level 2, Level 3, and Level 4 charging stations in the San Jose–Sunnyvale and Cupertino metropolitan areas, as shown in Fig.~\ref{fig:scenario1}. The peak EV demand occurs at Hour 19:38, when the total electrical load on the network is at its highest. At this time, the combined demand from Public Level 2 and DC Fast Charging stations reaches approximately 130 MW, representing the worst-case scenario for grid stress.}

\begin{figure}
    \centering 
\includegraphics[scale=0.36,trim=2cm 0.5cm 5cm 0.5cm,clip]{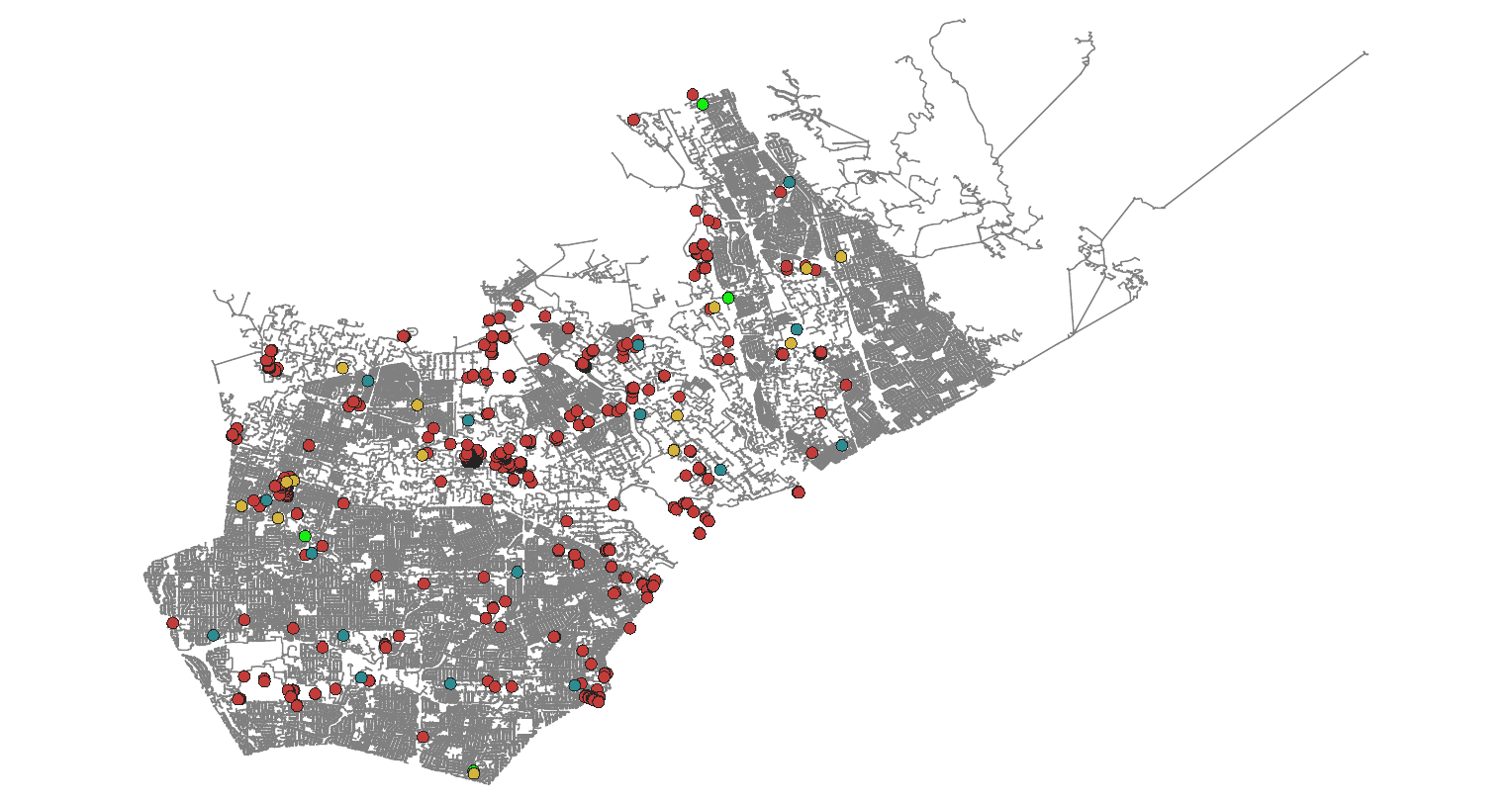}
\centering 
\caption{\mri{One line diagram of PU15 sub-region in San Jose, Sunnyvale,
and Cupertino. EV charging stations are determined as: Red Points: below 50 kW (Level 1), Yellow points:between 50 kW and 150 kW (Level 2),Blue points: between 150 kW and 350 kW (Level 3), Green Points: more than 350 kw (Level 4).}}
    \label{fig:addevtoPU15}
\end{figure}

\rn{The peak load is proportionally distributed across 951 EV charging stations based on their output capacities. Specifically, 895 stations operate below 50 kW, 24 stations fall within the 50–150 kW range, 18 stations operate between 150–350 kW, and 14 stations exceed 350 kW. Based on these classifications and their geographic locations, the assigned load per station is approximately 115.36 kW, 230.72 kW, 461.44 kW, and 922.9 kW for Level 1, 2, 3, and 4 stations, respectively. The load is proportionally distributed among the stations according to their rated capacities. The loads are automatically added to the system and shown in Fig.~\ref{fig:addevtoPU15}. The load assignment process is fully automated, and the underlying code is generalized to accommodate a wide range of scenario conditions.}

\rn{Using the assigned load demands, we apply Algorithm~\ref{alg:load_assignment} to analyze the system in terms of power flows, power losses, and voltage variations across the network. Before integrating the EV charging loads, the total system demand was approximately 1697 MW. After adding the EV loads, the total demand increased to 1827 MW, an increase of about 7.67\%. This highlights the additional stress placed on the distribution network due to EV integration. Regarding active power losses, the results show a notable increase after the EV loads are added. Initially, the system experienced losses of approximately 95,213.02 kW. After integration, these losses rose to 106,650 kW, marking an increase of about 12.05\%. This increase reflects higher energy dissipation in the network, indicating greater strain on system infrastructure and emphasizing the need for effective mitigation strategies.}

\begin{figure}
    \centering 
\includegraphics[scale=0.43,trim=0.1cm 0.1cm 0.1cm 0.1cm,clip]{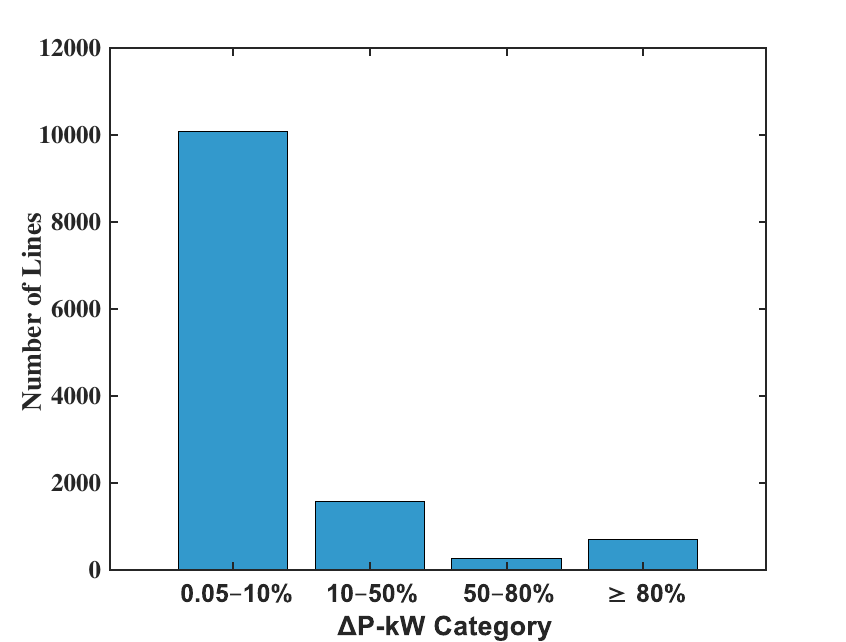}
\centering 
\caption{\rn{Distribution of power flow changes across distribution lines with 350,000 EVs integrated into the system. Most lines experience relatively small changes in power flow (0.05–10\%), while a smaller number of lines show moderate to high increases, highlighting potential stress points in the network.}}
    \label{fig:flow_hist350}
\end{figure}

\begin{figure}
    \centering 
\includegraphics[scale=0.43,trim=0.1cm 0.1cm 0.1cm 0.1cm,clip]{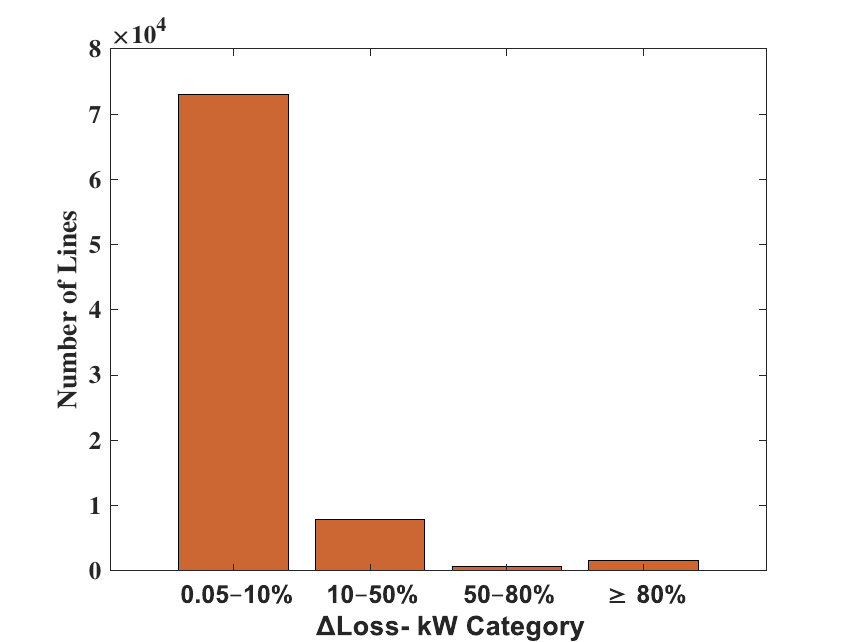}
\centering 
\caption{\rn{Distribution of power loss changes across distribution lines with the integration of 350,000 EVs. Most lines experience a modest increase in power loss (0.05–10\%), while fewer lines show moderate to significant increases, highlighting areas with higher stress due to EV-related loading.}}
    \label{fig:loss_hist350}
\end{figure}

\begin{figure}
    \centering 
\includegraphics[scale=0.4,trim=2cm 0.5cm 5cm 0.5cm,clip]{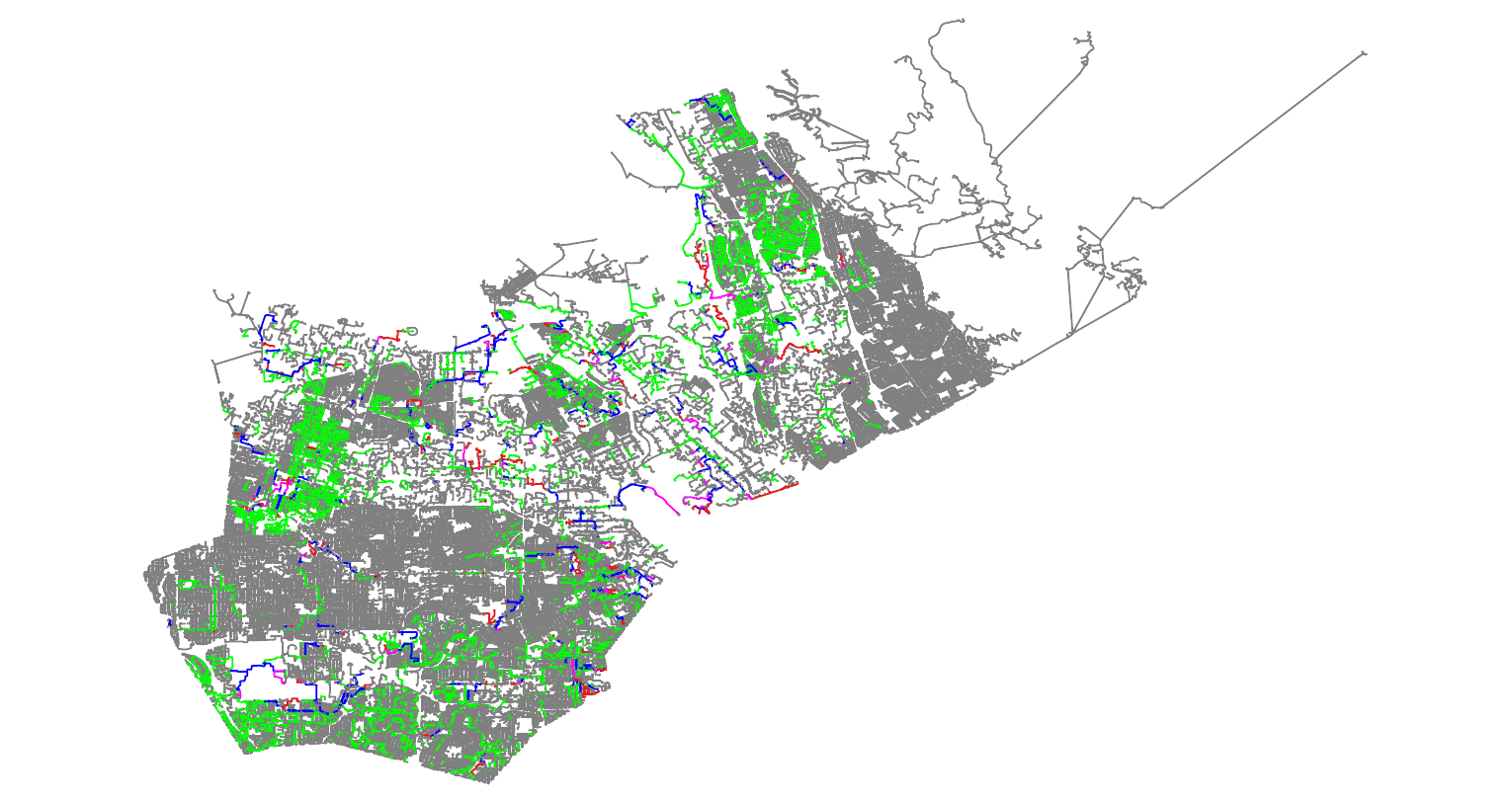}
\centering 
\caption{\rn{One-line diagram of the PU15 sub-region in San Jose, Sunnyvale, and Cupertino after integration of the load from 350,000 EVs. Line colors represent the percentage change in loading: gray for less than 0.05\%, green for 0.05-10\%, blue for 10-50\%, magenta for 50-80\%, and red for more than 80\%.}}
    \label{fig:add350evtoPU15}
\end{figure}

\rn{To evaluate the impact of EV integration on the distribution network, we analyzed the relative changes in power flows and power losses across the system. The histograms in Fig.\ref{fig:flow_hist350} and Fig.\ref{fig:loss_hist350} illustrate these changes by showing the distribution of percentage deviations between the original and EV-integrated scenarios. Each bin represents a range of normalized changes, and the height of each bin indicates how many lines experienced changes within that range. The power flow histogram shows that most lines experienced only moderate changes, while a smaller portion of the network saw significant increases in power flow due to the added EV charging demand. Specifically, the majority of lines (approximately 10,000) showed small changes in the range of 0.05–10\%, while several hundred to a few thousand lines experienced moderate (10–50\%) or high (50–80\%) increases. Around 1,000 lines exhibited changes greater than 80\%. This result is significant, as the proposed algorithm can easily identify the lines under stress. This is particularly useful for prioritizing grid upgrades and reinforcing critical components.
Similarly, the power loss histogram shows that some lines experienced noticeable increases in power losses, pointing to areas with higher thermal stress and reduced efficiency. These results help identify parts of the network most affected by EV integration and support decisions regarding grid upgrades and strategic placement of charging infrastructure. After assigning the EV station demand to the nearest network loads based on their GIS-based locations, we visualized the impact on the distribution network using a color-coded mapping scheme described in Algorithm~\ref{alg:QGIS}. The resulting map is shown in Fig.~\ref{fig:add350evtoPU15}.}

\begin{table*}[]
\caption{Key Factors and Assumptions for EV Charging Scenarios.}
\label{Table:FACTORS}
\begin{tabular}{|l|l|l|}
\hline
\textbf{Factor}& \textbf{States}& 
\textbf{Describtion}        \\ \hline
\begin{tabular}[c]{@{}l@{}}Plug-in \mri{EVs}\\  in the Fleet\end{tabular}&
Number  of \mri{EVs}& 
\begin{tabular}[c]{@{}l@{}}Based on national \mri{EV} stock projections from the U.S. Energy \\ Information Administration, the San Jose region could see 207,754 \mri{EV}\\ on the road by 2030. Alao there were approximately 61,230 plug-in \mri{EVs}\\  on the road   in the San Jose area as of the end of 2018.\end{tabular}          \\ \hline
\begin{tabular}[c]{@{}l@{}}Average Daily Miles \\ Traveled per vehicle\end{tabular} & 25 miles, 35 miles,  45 miles & \begin{tabular}[c]{@{}l@{}}Based on population density for this location, vehicles can be estimated to travel\\    25,35 or 45 miles per day on average.\end{tabular}    \\ \hline
\begin{tabular}[c]{@{}l@{}}Average Ambient \\ Temperature\end{tabular}& \begin{tabular}[c]{@{}l@{}}-4°F (-20°C), 14°F (-10°C), \\ 32°F (0°C), 50°F (10°C), \\ 68°F (20°C), 86°F (30°C)\\  104°F (40°C)\end{tabular}   & \begin{tabular}[c]{@{}l@{}}Based   on historical data for this location, a typical day averages 60.8°F, with\\  daily averages between 41°F and 80°F annually\end{tabular}       \\ \hline
\begin{tabular}[c]{@{}l@{}}Plug-in Vehicles that\\  are All-Electric\end{tabular}& 25\%,   50\% , 75\%         & \begin{tabular}[c]{@{}l@{}}This   filter has three options that let us choose a distribution of all-electric \\ vehicles (EVs), also called battery electric vehicles (BEVs). Each   distribution \\ is based on assumptions for the percentage of EVs and plug-in  hybrid electric \\ vehicles (PHEVs) that have enough energy storage to drive a   given distance.\\      25\% – PHEV Dominant:\\      This fleet has more PHEVs than EVs.\\      50\% – Equal Shares of EVs and PHEVs:\\      This fleet has about the same number of EVs as PHEVs.\\      75\% – EV Dominant.\\      This fleet has more EVs than PHEVs\end{tabular}     \\ \hline
\begin{tabular}[c]{@{}l@{}}Plug-in Vehicles that\\  are Sedans\end{tabular}         & 20\%, 50\%, 80\%          & \begin{tabular}[c]{@{}l@{}}This   filter has three options to choose the percentage of plug-in vehicles that\\  have a sedan body type.\end{tabular}           \\ \hline
\begin{tabular}[c]{@{}l@{}}Mix of Workplace \\ Charging\end{tabular}& \begin{tabular}[c]{@{}l@{}}20\%   Level 1 and 80\% Level 2 \\ 50\% Level 1 and 50\% Level 2\\  80\% Level 1 and 20\% Level   2\end{tabular} & \\ \hline
\begin{tabular}[c]{@{}l@{}}Access to Home \\ Charging\end{tabular}& \begin{tabular}[c]{@{}l@{}}50\%,   75\%, 100\% \\      with the following mix:\\      \\  20\% Level 1 and 80\% Level 2\\  50\%   Level 1 and 50\% Level 2 \\ 80\% Level 1 and 20\% Level\end{tabular} & \begin{tabular}[c]{@{}l@{}}This   filter has three options that let us choose the distribution of drivers in  \\  the fleet with access to charging at home:\\      50\% – 50\% have access to home charging and 50\% have no access\\      75\% – 75\% have access to home charging and 25\% have no access\\      100\% – 100\% have access to home charging and 0\% have no access\\   In addition, we can choose the distribution of charging power for drivers \\   with access\\       to home charging:\\      20\% Level 1 and 80\% Level 2\\      50\% Level 1 and 50\% Level 2\\      80\% Level 1 and 20\% Level 2\end{tabular}              \\ \hline
\begin{tabular}[c]{@{}l@{}}Preference for Home\\  Charging\end{tabular}       & 60\% 80\% 100\% & \begin{tabular}[c]{@{}l@{}}This   filter has three options that let us choose the distribution of drivers in   the fleet\\ who prefer home as their primary charging location with charging at   other locations\\  (including work) being used only as necessary to maintain  electric range:\\      \\      60\% – 60\% prefer primarily charging at home\\      80\% – 80\% prefer primarily charging at home\\      100\% – 100\% prefer primarily charging at home\end{tabular}           \\ \hline
\begin{tabular}[c]{@{}l@{}}Home Charging \\ Strategy\end{tabular}& \begin{tabular}[c]{@{}l@{}}Immediate – as fast as possible.\\  Immediate – as slow as possible.\\  Delayed – finish by departure. \\  Delayed – start at midnight\end{tabular} & \begin{tabular}[c]{@{}l@{}}This   filter has four options for home charging that let us explore load   flexibility by \\ adjusting the charging strategy.\\  *Immediate – as fast as possible:\\  This option assumes vehicles begin charging as soon as possible upon   arriving at \\ a charging location and charge at full power/speed until fully   charged or the vehicle \\ departs.\\  *Immediate – as slow as possible (even spread):\\  This option assumes vehicles begin charging immediately upon arriving at a   \\ charging location, but the charging speed/power is controlled to be as   slow/low \\ as possible to spread the charge evenly over the time the vehicle is   parked.\\  *Delayed – finish by departure:\\  This option assumes vehicles wait as long as possible to begin charging so   \\ they can still receive a full charge. This strategy uses arrival and   departure times\\  from the travel data referenced in the assumptions to shift   load during simulations.\\  *Delayed – start at midnight (only home charging):\\  This option assumes vehicles begin home charging at midnight because some  \\  vehicle owners elect to program their vehicles to start charging at a   specific time \\ overnight, which is often midnight.\end{tabular} \\ \hline
\begin{tabular}[c]{@{}l@{}}Workplace Charging \\ Strategy\end{tabular}              & \begin{tabular}[c]{@{}l@{}}Immediate   – as fast as possible.\\  Immediate – as slow as possible.\\ Delayed – finish by departure\end{tabular} & \begin{tabular}[c]{@{}l@{}}This   filter three options for workplace charging that let us explore load   flexibility\\  by adjusting the charging strategy same the home charging   strategy\end{tabular}       \\ \hline
\end{tabular}
\end{table*}

\subsection{Second Scenario}
\label{sec:second scenario}

\begin{figure}
    \centering 
\includegraphics[scale=0.45,trim=0.1cm 0.1cm 0.1cm 0.1cm,clip]{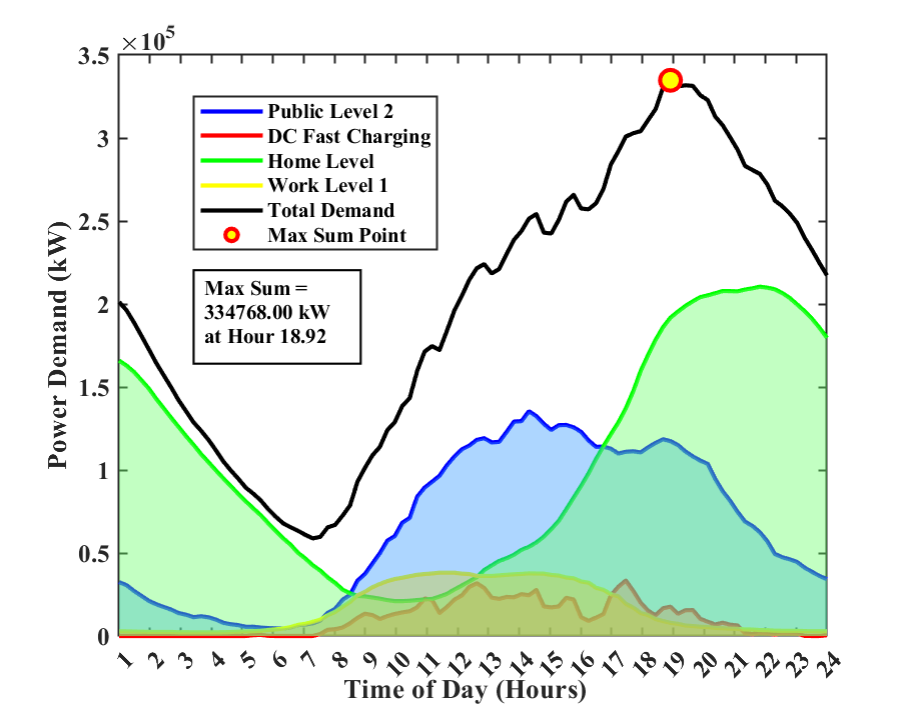}
\centering 
\caption{Charging demand profile over 24 hours for Scenario 2, with a projected fleet of 700,000 EVs in the target region. The load includes Level 1, Level 2, and DC Fast Charging (150 kW) stations. The peak demand reaches approximately 334.77 MW at 18:54, with average loads of 297 kW, 594 kW, 1188 kW, and 2376 kW assigned to Levels 1 through 4 stations, respectively.}
    \label{fig:seconddemand}
\end{figure}

\rn{In the second scenario, we increased the number of EVs to assess the impact of projected EV growth on the electric network in the target region by 2030. Specifically, we considered a fleet of 700,000 EVs. Based on this assumption, we generated load demand curves for Level 1, Level 2, and DC Fast Charging (150 kW) stations in the proposed areas, as shown in Fig.~\ref{fig:seconddemand}. In this case, the peak load demand occurs at 18:54 and reaches 334.77 MW. The assigned load per station is approximately 297 kW, 594 kW, 1188 kW, and 2376 kW for Levels 1, 2, 3, and 4 charging stations, respectively.}

\begin{figure}
    \centering 
\includegraphics[scale=0.43,trim=0.1cm 0.1cm 0.1cm 0.1cm,clip]{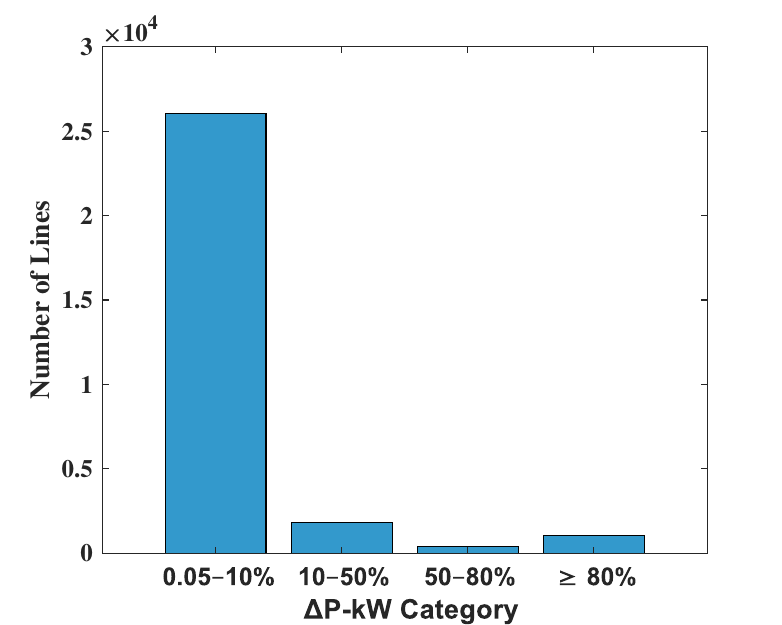}
\centering 
\caption{ Power flow change across distribution lines with 700,000 EVs. Compared to the first scenario, more lines experience increased power flow, with the majority falling in the 0.05–10\% change range. This highlights the growing impact on grid infrastructure under projected EV growth.}
    \label{fig:flow_hist}
    \vspace{-.4cm}
\end{figure}
\begin{figure}
    \centering 
\includegraphics[scale=0.43,trim=0.1cm 0.1cm 0.1cm 0.1cm,clip]{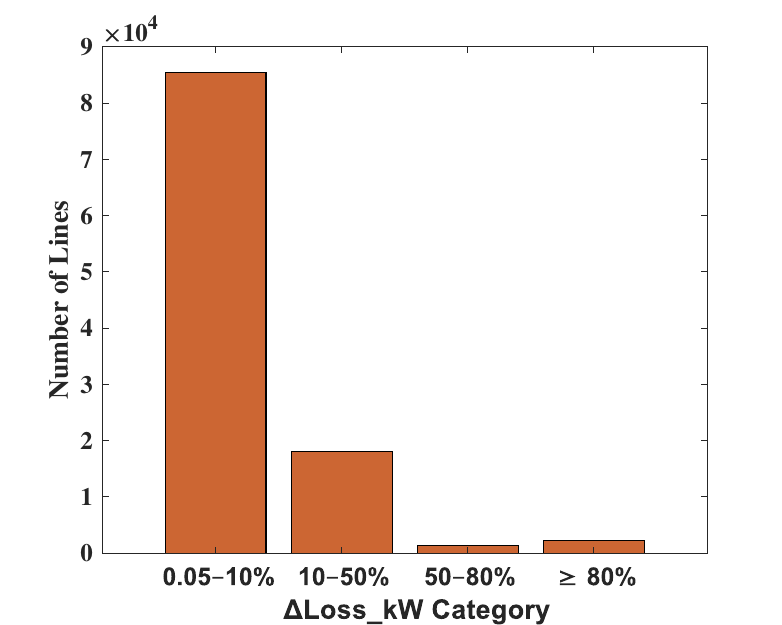}
\centering 
\caption{\rn{Distribution of power loss changes across distribution lines after integrating 700,000 EVs. The histogram categorizes the percentage increase in active power losses, showing that while most lines fall within the 0.05–10\% range, a significant number experience moderate to severe increases.}}
    \label{fig:loss_hist}
    \vspace{-.4cm}
\end{figure}
\begin{figure}
    \centering 
\includegraphics[scale=0.4,trim=2cm 0.5cm 5cm 0.5cm,clip]{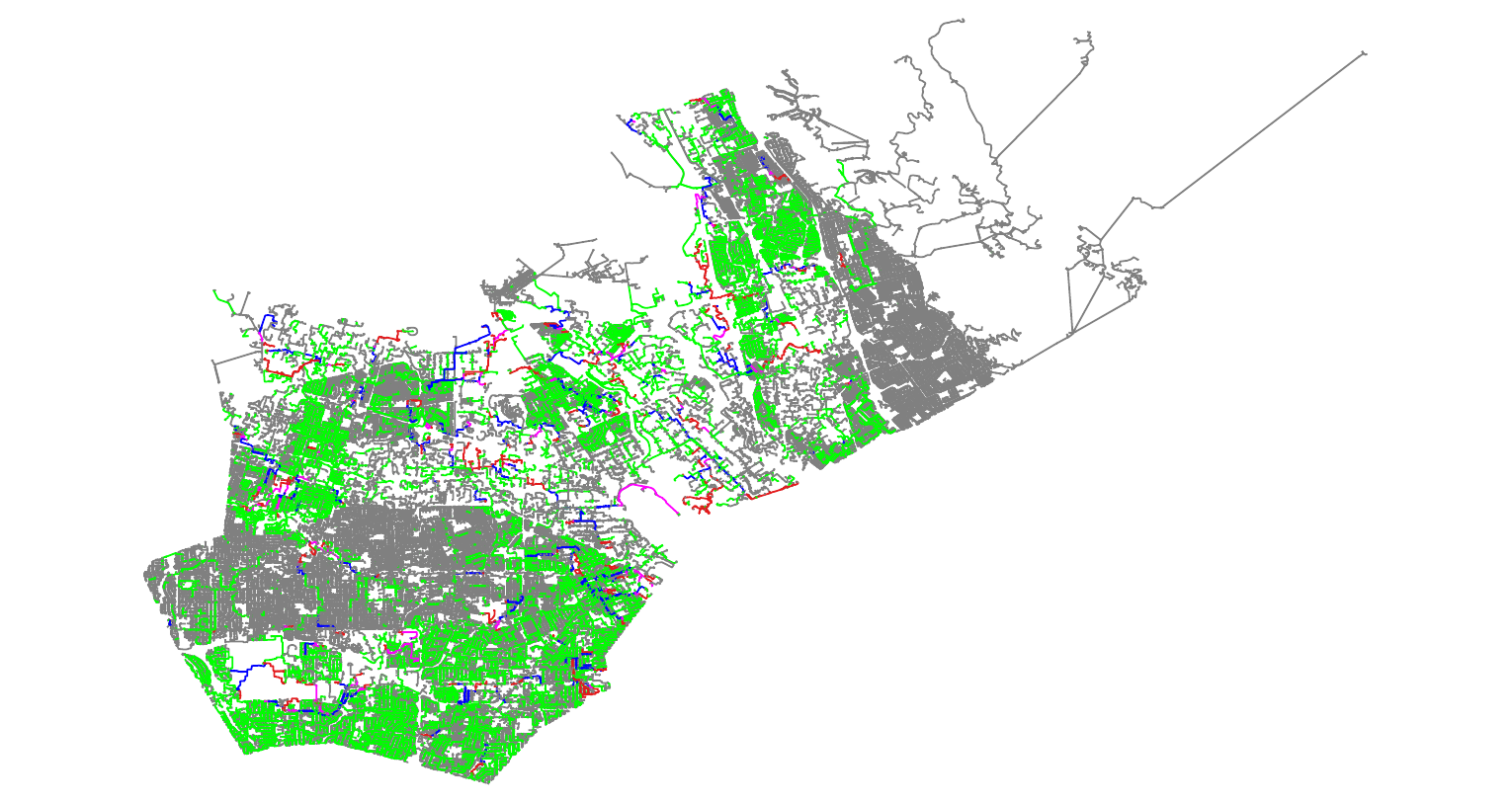}
\centering 
\caption{\rn{One-line diagram of the PU15 sub-region in San Jose, Sunnyvale, and Cupertino after integrating the load from 700,000 EVs. Line colors indicate the percentage change in loading: gray for <0.05\%, green for 0.05–10\%, blue for 10–50\%, magenta for 50–80\%, and red for $\geq $80\%.}}
    \label{fig:add700evtoPU15}
\end{figure}

\rn{Before integrating EV station loads, the total system demand was approximately 1697 MW. After adding the EV charging demand, the total system demand rose to 2032 MW, representing an increase of about 19.68\%. This highlights the additional burden placed on the distribution network due to the integration of EV charging infrastructure. Active power losses also increased significantly. Initially, the system experienced losses of approximately 95,213.02 kW. After incorporating the EV station loads, losses rose to 129,735.67 kW, an increase of about 36.32\%. This substantial rise indicates greater energy dissipation and increased stress on the distribution network, reinforcing the need for mitigation strategies. As in the previous scenario, the impact is visualized through histogram plots and a GIS-based map shown in Fig.\ref{fig:flow_hist}, Fig.\ref{fig:loss_hist}, and Fig.\ref{fig:add700evtoPU15}. In Fig.\ref{fig:flow_hist}, although most lines still fall in the 0.05–10\% change range, the number of lines in this bin increases from approximately 10,000 to about 26,000. The number of lines in the moderate-change range (10–50\%) roughly doubles, and the highest-impact category ($\geq$ 80\%) more than triples. Comparing Fig.\ref{fig:flow_hist} with Fig.\ref{fig:flow_hist350}, it is evident that doubling the number of EVs shifts more lines into the medium, and high-impact categories. While many lines still experience small changes, that group shrinks as more lines show substantial impacts.}

\rn{While most lines experience only modest changes, the extended tails in both the power flow and loss histograms indicate that higher EV penetration leads to a disproportionate increase in stress on a significant subset of feeders. Capturing these worst-case impacts, and how they scale with EV growth, requires high, resolution simulation and modeling. This level of detail is essential for identifying vulnerable lines and supporting targeted mitigation measures, such as network reconfiguration, or demand response strategies, to maintain system reliability under heavy EV loading.}

\rn{In the 350,000-EV scenario (Fig.\ref{fig:add350evtoPU15}), our framework identifies a limited number of distribution lines (highlighted in yellow) that show noticeable deviations in power flow compared to the base case (Fig.\ref{fig:addevtoPU15}). When EV penetration increases to 700,000 vehicles (Fig.~\ref{fig:add700evtoPU15}), the number and spatial extent of these high-change lines grow significantly. A larger portion of the network experiences variations in power flow, indicating that increased uncontrolled charging amplifies stress not only on individual circuits but also alters overall loading patterns across the grid. Although the current power distribution system can tolerate the present volume of EV loads, a substantial increase in EV adoption, such as that projected in the second scenario, can impose significant additional stress on the grid, potentially leading to congestion in distribution lines and overloading of transformers. This underscores the importance of proactive planning.}

\rn{By automatically generating GIS-based visualizations and quantifying per-line flow changes, the proposed approach enables planners to identify emerging stress points and evaluate mitigation strategies before costly upgrades or reliability issues arise. Moreover, the framework is easily adaptable to varying load demands, making it a valuable tool for both planning and operational decision-making in modern power systems.}


\section{Conclusion}
\label{sec:conclusion}

\rn{This paper presents a high-fidelity modeling framework to evaluate the impact of EVs on power distribution systems. Using publicly available data from the SMART-DS dataset, we developed a detailed and geographically grounded model of distribution feeders in the San Francisco region. Simulations are performed using the open-source tool OpenDSS, enabling accurate and scalable analysis of EV integration within the electric grid.
To demonstrate the capabilities of the model, we analyzed multiple scenarios, including both controlled and uncontrolled charging patterns, and varying levels of EV penetration. The results show that unmanaged charging can lead to significant overloading of grid components, while spreading out vehicle charging over time can substantially reduce stress on the system.
The proposed model is highly flexible and can be adapted to a wide range of demand profiles and infrastructure configurations. It is well-suited for both planning and operational decision-making in power systems. By identifying vulnerable components, such as overloaded lines and transformers, under projected EV load growth, the model supports proactive grid reinforcement and investment prioritization.
Additionally, this framework offers a valuable tool for researchers and utility planners to simulate future EV adoption scenarios, assess grid resilience, and evaluate the effectiveness of mitigation strategies such as optimized charging schedules and strategic placement of charging stations. As EV adoption continues to rise, such tools will be essential for ensuring reliable and efficient operation of the electric distribution network.}

\bibliographystyle{IEEEtran}
\IEEEtriggeratref{50}
\bibliography{ref}
\end{document}